\documentclass[10pt]{amsart}
\usepackage{amssymb,amsmath,amsfonts,epsfig,amscd}
\begin{document}
\bibliographystyle{plain}
\newtheorem{thm}{Theorem}[section]
\newtheorem{lem}[thm]{Lemma}
\newtheorem{prop}[thm]{Proposition}
\newtheorem{cor}[thm]{Corollary}
\newtheorem{conj}[thm]{Conjecture}
\newtheorem{mainlem}[thm]{Main Lemma}
\theoremstyle{definition}
\newtheorem{defn}[thm]{Definition}
\newtheorem{remark}[thm]{Remark}

\newcommand{\R}{\mathbb R}
\newcommand{\C}{\mathbb C}
\newcommand{\Q}{\mathbb Q}
\newcommand{\N}{\mathbb N}
\newcommand{\Z}{\mathbb Z}
\newcommand{\Nc}{\mathcal N}

\newcommand{\PSL}{\operatorname{PSL}}
\newcommand{\Sl}{\operatorname{Sl}}
\newcommand{\Core}{\operatorname{Core}}
\newcommand{\arrow}{\rightarrow}
\newcommand{\defeq}{\stackrel{\text{def}}{=}}

\newtheorem*{namedtheorem}{\theoremname}
\newcommand{\theoremname}{testing}
\newenvironment{named}[1]{\renewcommand{\theoremname}{#1}\begin{namedtheorem}}{\end{namedtheorem}}
\def\square{\hfill${\vcenter{\vbox{\hrule height.4pt \hbox{\vrule width.4pt height7pt \kern7pt \vrule width.4pt} \hrule height.4pt}}}$}

\newenvironment{pf}{{\it Proof:}\quad}{\square \vskip 12pt}
\newcommand{\U}{\ensuremath{\widetilde}}
\newcommand{\Hn}{\ensuremath{\mathbb{H}^3}}
\newcommand{\trivgp}{\text{id}}
\title{Finitely generated subgroups of lattices in $\PSL_2 \mathbb{C}$}
\author{Yair Glasner, Juan Souto, and Peter Storm}
\date{April 7th, 2005}
\begin{abstract}
Let $\Gamma$ be a lattice in $\PSL_2 (\C)$. The pro-normal topology on $\Gamma$ is defined by
taking all cosets of non-trivial normal subgroups as a basis. This topology is finer than the
pro-finite topology, but it is not discrete. We prove that every finitely generated subgroup
$\Delta < \Gamma$ is closed in the pro-normal topology. As a corollary we deduce that if $H$
is a maximal subgroup of a lattice in $\PSL_2( \mathbb{C})$ then either $H$ is finite index or
$H$ is not finitely generated.
\end{abstract}
\thanks{The third author was partially supported by a National Science Foundation Postdoctoral Fellowship.}

\maketitle
\section{Introduction}
In most infinite groups, it is virtually impossible to understand the lattice of all
subgroups. Group theorists therefore focus their attention on special families of subgroups
such as finite index subgroups, normal subgroups, or finitely generated subgroups. This paper
will study the family of finitely generated subgroups of lattices in $\PSL_2(\C)$.  Recall
that lattices in $\PSL_2 (\C)$ are the fundamental groups of finite volume hyperbolic
$3$-orbifolds.  For these groups we establish a connection between finitely generated
subgroups and normal subgroups.  This connection is best expressed in the following
topological terms.

If a family of subgroups $\mathcal{N}$ is invariant under conjugation and satisfies the condition,
\begin{equation}
\label{cond:well_deff} \text{for all } N_1,N_2 \in \mathcal{N} \ \ \exists \ N_3 \in
\mathcal{N} {\text{ such that }} N_3 \le N_1 \cap N_2
\end{equation}
then one can define an invariant topology in which the given family of groups constitutes a basis of open
neighborhoods for the identity element.
The most famous example is the {\it
pro-finite} topology, which is obtained by taking $\mathcal{N}$ to be the family of (normal)
finite index subgroups. When the family of all non-trivial normal subgroups satisfies
condition (\ref{cond:well_deff}), we refer to the resulting topology as the {\it pro-normal}
topology. The pro-normal topology is usually much finer than the pro-finite topology.

With this terminology we can state the main theorem.
\begin{thm} \label{thm:main}
Let $\Gamma$ be a lattice in $\PSL_2(\C)$.  Then the pro-normal topology is well defined on
$\Gamma$ and every finitely generated subgroup $\Delta < \Gamma$ is closed in this topology.
Moreover if $\Delta$ is of infinite index then it is the intersection of open subgroups
strictly containing $\Delta$.
\end{thm}

{\noindent}Recall that a subgroup is open in the pro-normal topology if and only if it contains a nontrivial normal subgroup, an open subgroup is automatically closed, and a group is closed in the pro-normal topology if and only if it is the intersection of open subgroups (see Remark \ref{rem:closed_subgroups}).  From this it follows that the second statement of Theorem \ref{thm:main} is stronger only if $\Delta$ is open.  As we will see in the proof, ``most'' finitely generated subgroups are not open.  Nevertheless, the stronger
version is used below to prove quickly Corollary \ref{cor:max}.

It is a well known question whether or not Theorem \ref{thm:main} remains true for the coarser
pro-finite topology. Groups in which every finitely generated group is closed in the
pro-finite topology are called {\it locally extended residually finite}, or {\it LERF} for
short. The list of groups that are known to be LERF is rather short. It includes lattices in
$\PSL_2(\R)$ (see \cite{Scott:LERF}) and a few examples of lattices in $\PSL_2(\C)$ (see
\cite{Gitik:3mfld_LERF}).

Let $\Gamma$ be a group $\Delta < \Gamma$ a subgroup. We say that $\Delta$ is a {\it maximal}
subgroup if there is no subgroup $\Sigma$ such that $\Delta \lneqq \Sigma \lneqq \Gamma$. In
\cite{MS:Maximal}, Margulis and So\u{\i}fer prove that every finitely generated linear group
that is not virtually solvable admits a maximal subgroup of infinite index.
On various occasions both Margulis and So\u{\i}fer asked the question, in which cases can such
maximal subgroups of infinite index be finitely generated? Easy examples of finitely generated
maximal subgroups such as $\PSL_2(\Z) < \PSL_2(\Z[1/p])$ indicate that one should be careful
about the exact phrasing of the question. So\u{\i}fer suggested the setting of lattices in
simple Lie groups, and indeed even for $\PSL_3(\Z)$ the question is wide open.

Let $\Gamma$ be a group $\Delta < \Gamma$ a subgroup. We say that $\Delta$ is {\it pro-dense}
if $\Delta$ is a proper subgroup and $\Delta N = \Gamma$ for every normal subgroup $\trivgp
\ne N \lhd \Gamma$. When the pro-normal topology on $\Gamma$ is well defined pro-dense
subgroups are exactly the subgroups that are dense in this topology. In \cite{GG:Primitive} it
is proven that if $\Gamma$ is a finitely generated linear group with simple Zariski closure,
or if it is a non-elementary hyperbolic group with no finite normal subgroups, then $\Gamma$
admits a pro-dense subgroup.  Note that all lattices in $\PSL_2 (\C)$ satisfy the first
hypothesis, and uniform lattices satisfy both hypotheses.

In groups that have few quotients such as simple groups it can be relatively easy to construct
pro-dense subgroups. This becomes more of a challenge in groups that have many normal
subgroups, such as hyperbolic groups. It is conjectured in \cite{GG:Primitive} that a
pro-dense subgroup of a hyperbolic group cannot be finitely generated.

Our main theorem above answers these two questions for the special
case of lattices in the group $\PSL_2(\C)$.

\begin{cor} \label{cor:max}
Let $\Gamma$ be a lattice in $\PSL_2(\C)$ and $\Delta < \Gamma$ a maximal subgroup of infinite
index or a pro-dense subgroup. Then $\Delta$ cannot be finitely generated.
\end{cor}
\begin{proof}
In both cases the only open subgroup strictly containing $\Delta$ is $\Gamma$ itself. If
$\Delta$ is maximal this is true by definition, even without the openness condition. If
$\Delta$ is pro-dense then it is dense in the pro-normal topology and it cannot be contained
in an open (and hence closed) proper subgroup. Thus $\Delta$ cannot be finitely generated or
it would violate the conclusion of Theorem \ref{thm:main}.
\end{proof}

The main new tool used to prove Theorem \ref{thm:main} and Corollary \ref{cor:max} is the recent resolution of the Marden conjecture \cite{Agol:Tameness, CG:Tameness}.  This powerful result is very specific to hyperbolic $3$-space.  It therefore seems likely that new techniques will be required to tackle these questions for other simple Lie groups.

The authors would like to thank I. Agol, I. Kapovich, G. Margulis, P. Schupp, and G.
So\u{\i}fer.

\section{Geometrically finite groups.}
\label{sec:GF}
\begin{prop} \label{main prop}
Let $\Gamma < \text{Isom}^+ (\mathbb{H}^m)$ be a lattice (for $m>2$). Let $\Delta < \Gamma$ be
a geometrically finite subgroup of infinite index.  Then there exists a normal subgroup
$\trivgp \ne N \lhd \Gamma$ such that $N \cap \Delta = \trivgp$.
\end{prop}
By Selberg's lemma, there exists a finite index torsion free normal subgroup
$\Gamma^{\text{tf}} \le \Gamma$.  Let $M$ be the hyperbolic manifold $\mathbb{H}^m /
\Gamma^{\text{tf}}$.  For a hyperbolic element $g \in \Gamma^{\text{tf}}$, let $g^* \subset M$
denote the closed geodesic corresponding to the conjugacy class of $g$.  Let $A_g \subset
\mathbb{H}^m$ denote the axis of the isometry $g$.  Define the torsion free subgroup
$\Delta^{\text{tf}} := \Delta \cap \Gamma^{\text{tf}}$.

To handle lattices with torsion, it will be necessary to consider the lattice $\text{Aut} :=
\text{Aut}(\Gamma^{\text{tf}}) < \text{Isom}^+ (\mathbb{H}^m)$, which by Mostow rigidity is a
finite extension of $\Gamma$, and the finite group $\text{Out} :=
\text{Out}(\Gamma^{\text{tf}})$ of outer-automorphisms of $\Gamma^{\text{tf}}$ acting
isometrically on $M$.  The reader is encouraged to assume on a first reading that $\Gamma$ is
torsion free and $\text{Out}(\Gamma)$ is the trivial group.

\begin{lem} \label{lemma (1)}
There exists a primitive hyperbolic element $\gamma \in \Gamma^{\text{tf}}$ whose axis
$A_\gamma$ cannot be translated into the convex hull of $\Delta^{\text{tf}}$ by an element of
$\text{Aut}$. Moreover, there exist constants $\eta, L>0$ such that: for any element $\phi \in
\text{Aut}$ and any hyperbolic element $\delta \in \Delta^{\text{tf}}$ the intersection
$$\mathcal{N}_\eta A_\delta \cap A_{\phi. \gamma},$$ where $\mathcal{N}_\eta A_{\delta}$ is a
radius $\eta$ neighborhood of $A_{\delta}$, is a (possibly empty) geodesic segment of length
less than $L$.
\end{lem}
\begin{pf}
For a Riemannian manifold $V$ let $SV$ denote the unit tangent bundle of $V$.  Let $\U{M}$
denote the hyperbolic manifold $\mathbb{H}^m / \Delta^{\text{tf}}$, and let $\iota: S\U{M}
\longrightarrow SM$ denote the canonical covering map.

Let us first assume that $\Gamma^{\text{tf}} = \text{Aut}$.  As a first step, under this
assumption we will find a hyperbolic element $\gamma \in \Gamma^{\text{tf}}$ such that no lift
of $\gamma^* \subset M$ to $\mathbb{H}^m$ is contained in the convex hull of
$\Delta^{\text{tf}}$.

Let $\Lambda \subset \partial \mathbb{H}^m = \mathbb{S}^{n-1}$ be the limit set of
$\Delta^{\text{tf}}$.  The exponential map provides a natural homeomorphism $\mathbb{H}^m
\times \mathbb{S}^{n-1} \longrightarrow S\mathbb{H}^m$.  Let $\U{X} \subset S \mathbb{H}^m$ be
the image of $\text{Hull}(\Delta^{\text{tf}}) \times \Lambda$.  $\U{X}$ is a closed subset of
measure zero because $\Delta^{\text{tf}}$ is geometrically finite.  $\U{X}$ is
$\Delta^{\text{tf}}$-invariant and thus projects to a closed measure zero subset $X$ of $S
\U{M}$.

Let $C \subset \U{M}$ be the convex core of $\U{M}$.  $X$ lives ``over'' $C$.  Pick a small
closed ball $U \subset SM$ such that $\iota^{-1} (U) \cap SC$ is contained in the interior of
$SC$.  $U$ has only a finite number of preimages intersecting $SC$.  Therefore $\iota(
\iota^{-1} (U) \cap X) \subset SM$ is a closed measure zero set.  Closed geodesics form a
dense subset of $SM$.  So we may pick a $\gamma \in \Gamma^{\text{tf}}$ such that $\gamma^*
\subset SM$ (lift $\gamma^*$ to $SM$ in the natural way) intersects $U$ and is disjoint from
$\iota( \iota^{-1} (U) \cap X)$.  This precisely says that no lift of $\gamma^*$ to
$\mathbb{H}^m$ will be contained in the convex hull of $\Delta^{\text{tf}}$.  The completes
the first step of the proof under the assumption that $\Gamma^{\text{tf}} =
\text{Aut}(\Gamma^{\text{tf}})$.

In the case where $\text{Aut}$ properly contains $\Gamma^{\text{tf}}$, we must find a
hyperbolic element $\gamma \in \Gamma^{\text{tf}}$ such that for any element $\phi \in
\text{Out}$, no lift of the closed geodesic $\phi(\gamma^*) \subset M$ is contained in the
convex hull of $\Delta^{\text{tf}}$.  Since $\text{Out}$ is finite, and a finite union of
closed measure zero sets is closed and measure zero, this can be proved by repeating the above
argument a finite number of times to the subgroups $\{ \phi. \Delta^{\text{tf}} \}_{\phi \in
\text{Out}}$.  This completes the first step of the proof.

Fix the element $\gamma \in \Gamma^{\text{tf}}$ from the first step.  Let $r > 0$ be the
infimal injectivity radius of $M$ over the points of the closed geodesic $\gamma^* \subset M$.

The second step is to prove the constants $\eta, L >0$ exist for the element $\gamma$. Suppose
these constants do not exist.  Then there exist sequences $\eta_i \searrow 0$ (assume $\eta_i
< 1$), $L_i \rightarrow \infty$, elements $\phi \in \text{Aut}$, and hyperbolic elements
$\delta_i \in \Delta^{\text{tf}}$ such that
$$\mathcal{N}_{\eta_i} A_{\delta_i} \cap A_{\phi. \gamma}$$
is a geodesic segment $\sigma_i$ of length $L_i$.  Let $\mathcal{F} \subset \mathbb{H}^m$ be
the compact subset of a fundamental domain for $\Delta^{\text{tf}}$ given by points within
distance $1$ from the convex hull of $\Delta^{\text{tf}}$, which cover a point in $\U{M}$ of
injectivity radius at least $r/2$.  For each $i$ there is an $h_i \in \Delta^{\text{tf}}$ such
that $h_i$ moves the midpoint of the geodesic segment $\sigma_i$ into $\mathcal{F}$.  Since
the group $\text{Aut}$ is discrete and $\mathcal{F} \subset \mathbb{H}^3$ is compact, we may
pass to a subsequence where the segments $h_i (\sigma_i)$ all lie on a common bi-infinite
geodesic $\sigma$.  Since $\eta_i \rightarrow 0$, the geodesic $\sigma$ must lie in the convex
core of $\Delta^{\text{tf}}$.  But $\sigma$ is the axis of a conjugate of $\gamma$ (in the
group $\text{Aut}$).  This contradicts the first part of the proof, and completes the proof of
the lemma.
\end{pf}

Now fix a hyperbolic element $\gamma \in \Gamma^{\text{tf}}$ satisfying Lemma \ref{lemma (1)}.
Let $\ell>0$ be the length of the closed geodesic $\gamma^*$.  Let $\text{Out}.\gamma^*$
denote the immersed $1$-manifold $\cup_{\phi \in \text{Out}} \phi(\gamma^*)$.  An annoying
detail is that $\text{Out}.\gamma^* \subset M$ may not be an embedded $1$-manifold.
Nonetheless, there are constants $\rho>0$, $c \in (0,1)$, and a finite union of embedded
intervals $\{ I_\alpha \}_\alpha \subseteq \gamma^* \subset M$ of total length at least $c
\cdot \ell$ such that: a $\rho$-tubular neighborhood of $I_\alpha$ is embedded and disjoint
from a $\rho$-tubular neighborhood of $\phi(I_\beta)$ whenever $\phi(I_\beta) \neq I_\alpha$
as subsets of $M$.

Define $N_n \le \Gamma^{\text{tf}}$ to be the subgroup generated by the set
$$\bigcup_{\phi \in \text{Aut}} \phi.\gamma^n.$$
Notice that $N_n$ is a characteristic subgroup of $\Gamma^{\text{tf}}$. Since $\Gamma^{\text{tf}}$ is a normal subgroup of $\Gamma$, $N_n$ is also a normal subgroup of $\Gamma$.  This subgroup $N_n$ has the advantage of being torsion free. 

\begin{remark} \label{rem:characteristic}
In fact we could have initially chosen $\Gamma^{\text{tf}}$ to be a characteristic subgroup of $\Gamma$.  This is because every finite index subgroup in a finitely generated group contains a further finite index subgroup that is characteristic; just take the intersection of all subgroups of the given index.  In this way we can arrange for $N_n$ to be a characteristic subgroup of $\Gamma$.
\end{remark}

We are now prepared to prove Proposition \ref{main prop}.  The argument is similar to
\cite[5.5F]{Gromov_book}.   We will show that for a sufficiently large number $n$, the group
$N = N_n$ has a trivial intersection with $\Delta$.

Recall that $\ell$ is the length of the closed geodesic $\gamma^* \subset M$.  Recall the
constants $\eta$ and $L$ from Lemma \ref{lemma (1)}.  Define $\varepsilon := \text{min} \{
\rho, \eta \}$.  Fix an integer $n_0$ satisfying
$$c\ell \cdot n_0 > 2\pi/\varepsilon + L.$$
Pick any $n> n_0$.  Suppose that $\Delta \cap N_n$ is nontrivial.  Since $N_n \le
\Gamma^{\text{tf}}$ is torsion free and $\Delta^{\text{tf}} \le \Delta$ is a finite index
subgroup, it follows that there is a nontrivial element $\delta \in \Delta^{\text{tf}} \cap
N_n$.

For clarity let us break the proof into two cases, depending on whether or not $\Gamma$ has
parabolic elements.

\vskip 4pt

{\noindent}\emph{Case (1):}  Assume that $\Gamma$ has no parabolic elements. \vskip 2pt
Topologically the existence of an element $\delta \in \Delta^{\text{tf}} \cap N_n \le
\Gamma^{\text{tf}}$ says there is a compact genus zero surface $S$ with $(k+1)$ circular
boundary components $C_0, C_1, \ldots, C_k$ and a map
$$f: (S, C_0, C_1, \ldots, C_k) \longrightarrow (M, \delta^*, \phi_1(\gamma^*), \ldots, \phi_k(\gamma^*))$$
such that $\phi_j \in \text{Out}$, $f: C_0 \longrightarrow \delta^*$ is a parametrization of
$\delta^*$ and $f: C_j \longrightarrow \phi_j(\gamma^*)$ ($j>0$) wraps $i_j\cdot n$ times
around $\phi_j (\gamma^*)$ (possibly in the reverse direction) for some $i_j \in \mathbb{N}$.
Let us assume we have found such a surface $S$ with the minimal number of boundary components.

\begin{figure} \label{the figure}
\begin{center}
\epsfig{file=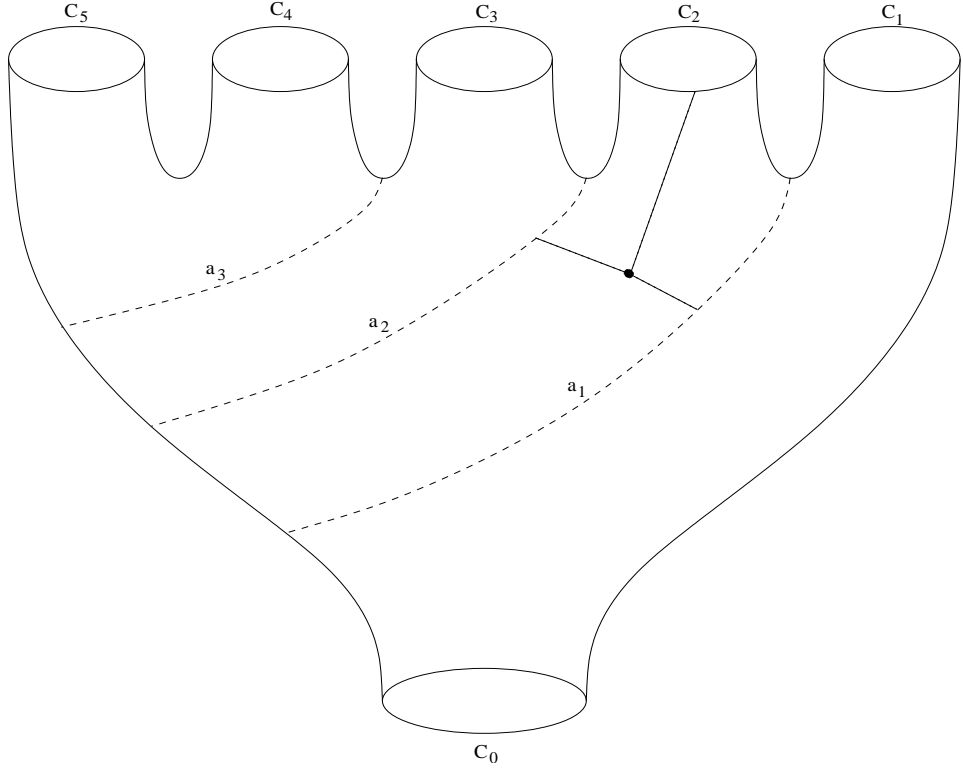,height=4in, width=5in} \caption{The surface S}
\end{center}
\end{figure}

We would now like to alter $f$ by a homotopy rel boundary to form a pleated surface (or a
minimal surface).  Consider Figure 1.  Add simple closed curves $\{ a_1,
a_2,\ldots, a_{k-2} \}$ to $S$ as shown in the figure, forming a pants decomposition of $S$.
If for some $j$ the curve $f(a_j)$ were homotopically trivial then we could cut along $a_j$,
add a disk, and find a surface with fewer boundary components than $S$.  Therefore by
minimality, each curve $f(a_j) \subset M$ is homotopic to a closed geodesic $a_j^*$.  By
another minimality argument one can show that for $1 \le j \le k-1$ the curve $f(a_j)$ is not
freely homotopic into any of the closed geodesics $\{ \phi (\gamma^*) \}_{\phi \in
\text{Out}}$.

For $1 \le j \le k-3$ consider the pants bounded by $a_j, a_{j+1},$ and $C_{j+1}$.  Add a
basepoint
to the $j^{th}$ pair of pants $P_j$ and paths to the boundary components as shown in the
figure for $j=1$. These (together with a choice of orientation) make each boundary component
$a_j, a_{j+1},$ and $C_{j+1}$ a well defined element $[a_j], [a_{j+1}],$ and $[ C_{j+1}]$
respectively of the fundamental group.  Suppose $f_*([a_j])$ and $f_* ([a_{j+1}])$ are
contained in a common cyclic subgroup of $\pi_1 (M)$.  Then $f_* ([C_{j+1}])$ would be
contained in the same cyclic subgroup.  In this case $f(a_j)$ would be freely homotopic into
$\phi_{j+1}(\gamma^*)$, which was ruled out in the previous paragraph.  Therefore $f_*([a_j])$
and $f_* ([a_{j+1}])$ are not contained in a common cyclic subgroup of $\pi_1 (M)$.  A similar
minimality argument shows that any pair of the triple $a_{k-2}$, $C_{k-1}$, and $C_k$ is not
mapped by $f$ into a common cyclic subgroup of $\pi_1 (M)$.

With this information we may ``pull $f$ tight'' along the curves $a_1, a_2, \ldots, a_{k-2}$
to the closed geodesics $a_1^*, a_2^*, \ldots, a_{k-2}^* \subset M$, and then use ideal
hyperbolic triangles to map each pair of pants into $M$ as a pleated surface.  (See Thurston's
description of this in \cite[Sec.2]{Th3}.  Alternatively one can map each pair of pants to a
minimal surface with geodesic boundary.)  Let us denote the resulting map again by $f$. There
is a hyperbolic metric on $S$ with geodesic boundary components such that $f$ maps paths in
$S$ to paths in $M$ of the same length.  Let us now use $S$ to denote the surface equipped
with this hyperbolic metric.

Recall that on the embedded intervals $\{ I_\alpha \}_\alpha \subset \gamma^*$ of total length
at least $c \ell$, the $\rho$-tubular neighborhoods of the intervals $\{ \phi(I_\alpha)
\}_{\phi \in \text{Out}}$ are pairwise disjoint and embedded.  This implies that at least
$100c$ percent of $C_1$ (i.e. curves of total length at least $c i_j n \ell$) is at least a
distance $2\rho$ from $\cup_{i>1} C_i$.

For $0 \le i < k$ let $\sigma_i \subset S$ denote the geodesic segment in $S$ of shortest
length joining $C_i$ to $C_k$.  Note that by minimality, the $\sigma_i$ are pairwise disjoint.
Cut $S$ along the segments $\{ \sigma_i \}_{0 \le i \le k}$ to form the simply connected
surface $\text{cut}S$.  Lift $f$ to a map
$$\U{f} : \text{cut}S \longrightarrow \mathbb{H}^m.$$
Let $p \in C_1$ be the point closest to $C_0$.  By Lemma \ref{lemma (1)}, the (possibly
disconnected) geodesic segment $\U{f}(C_1 - B_S(p, L))$ is entirely outside an
$\varepsilon$-neighborhood of $\U{f}(C_0)$. Therefore a $\varepsilon$-neighborhood of the
geodesic segment $C_1 - B_S (p,L)$ is disjoint from $C_0$.

Perform the above procedure for each boundary component $C_i$ ($2\le i \le k$), where the role
of $C_k$ above is replaced with $C_1$.  This shows that the hyperbolic area of an
$\varepsilon$-neighborhood of $\cup_{i>1} C_i$ is at least
$$(nc\ell - L) \cdot \varepsilon \cdot k.$$
By Gauss-Bonnet, the area of $S$ is $-2\pi (2-k-1)=2\pi k - 2 \pi.$  (In the case of a minimal
surface the area is at most this number.)  Therefore
$$2\pi k  \ge nc\ell \varepsilon k - L \varepsilon k$$
yielding
$$n c\ell \le 2\pi/\varepsilon + L.$$
This contradicts our choice of $n> n_0$.  This proves (2) under the assumption that $\Gamma$
contains no parabolic elements. \vskip 4pt

{\noindent}\emph{Case (2):}  $\Gamma$ contains parabolic elements. \vskip 2pt If $\delta$ is a
parabolic isometry, then $S$ is replaced with a genus zero surface with $k$ boundary
components and one puncture.  The map $f$ then maps the end of $S$ out the cusp of $M$
corresponding to $\delta$.

Proceeding as in Case (1), add the curves $a_j$ to $S$.  It may happen that a curve $f(a_j)$
represents a parabolic element of $\Gamma$.  In this case cut $S$ along the first such curve
$a_j$ and map the resulting surface into $M$ taking the curve $f(a_j)$ out the corresponding
cusp of $M$.  Proceed with the above argument treating the new boundary curve $a_j$ as the
``last'' boundary curve $C_k$ above.

The argument of Case (1) goes through with obvious alterations on the resulting finite area
pleated surface, and the proof is left to the reader. \square \vskip 6pt

\section{Finitely generated subgroups} \label{sec:pf}
We begin by quoting the two main theorems from hyperbolic geometry which we will use.
\begin{thm} \label{Marden conjecture} (Marden Conjecture) \textnormal{[}Agol \cite{Agol:Tameness}, Calegari-Gabai \cite{CG:Tameness}\textnormal{]}
Let $\Delta$ be a discrete finitely generated torsion-free subgroup of $\PSL_2 \mathbb{C}$.
Then the hyperbolic manifold $\mathbb{H}^3 / \Delta$ is topologically tame, i.e. it is
homeomorphic to the interior of a compact manifold with boundary.
\end{thm}

\begin{thm} \label{Covering theorem} (Covering Theorem) \textnormal{[}Thurston \cite[Ch.9]{Th}/Canary \cite{Canary_covering}\textnormal{]}
Let $\Gamma < \PSL_2 \mathbb{C}$ be a torsion-free lattice.  Let $\Delta < \Gamma$ be a subgroup such that $\mathbb{H}^3/\Delta$ is topologically tame.  Then either \\
(1)  $\Delta$ is geometrically finite, or \\
(2)  there exists a finite index subgroup $H \le \Gamma$ and a subgroup $F \le \Delta$ of
index one or two such that $H$ splits as a semidirect product $F \rtimes \mathbb{Z}$.
\end{thm}

We may now combine Theorems \ref{Marden conjecture} and \ref{Covering theorem} with Selberg's
lemma to prove a strengthened version of the Covering Theorem for lattices with torsion.
\begin{cor} \label{covering cor}
Let $\Gamma < \PSL_2(\mathbb{C})$ be a lattice.  Let $\Delta < \Gamma$ be a finitely generated
subgroup.
Then either \\
(1)  $\Delta$ is geometrically finite, or \\
(2)  there exists a torsion free finite index subgroup $H \le \Gamma$ and a finite index
normal subgroup $F \unlhd \Delta$ such that $H$ splits as a semidirect product $F \rtimes
\mathbb{Z}$.
\end{cor}

\begin{pf}
Assume $\Delta$ is not geometrically finite.  Apply Selberg's lemma to obtain a torsion free
finite index normal subgroup $\Gamma^{\text{tf}} \le \Gamma$.  The subgroup $\Delta \cap
\Gamma^{\text{tf}}$ is normal torsion free and finite index in $\Delta$.  Apply Theorems
\ref{Marden conjecture} and \ref{Covering theorem} to the pair $\Gamma^{\text{tf}}$ and
$\Delta \cap \Gamma^{\text{tf}}$.

We know $\Delta \cap \Gamma^{\text{tf}}$ is not geometrically finite.  So there is a finite
index subgroup $H \le \Gamma^{\text{tf}}$ and a subgroup $F \le \Delta \cap
\Gamma^{\text{tf}}$ of index one or two such that $H$ splits as a semidirect product $F
\rtimes \mathbb{Z}$.  Since its index is at most two, $F$ is in fact a characteristic subgroup
of $\Delta \cap \Gamma^{\text{tf}}$, implying it is a normal subgroup of $\Delta$.  This
proves the corollary.
\end{pf}

We are now ready to prove our main theorem, which we restate here for convenience.
\begin{named}{Theorem \ref{thm:main}} \itshape
Let $\Gamma$ be a lattice in $\PSL_2(\C)$ then, the pro-normal topology is well defined on
$\Gamma$ and every finitely generated subgroup $\Delta < \Gamma$ is closed in this topology.
Moreover if $\Delta$ is of infinite index then it is the intersection of open subgroups
strictly containing $\Delta$.
 \normalfont
 \end{named}

 \begin{remark} \label{rem:closed_subgroups} Before beginning the proof let us recall why a subgroup is closed in the pro-normal topology if and only if it is the intersection of open subgroups. One direction is obvious, an open subgroup is also closed because its complement is a union of open cosets. Conversely, applying the definition of the topology, a set $S \subseteq \Gamma$ is open if and only if
 $$S = \bigcup_{N \lhd \Gamma} \ \bigcup_{ \{ \gamma \, | \, \gamma N \subset S \}} \gamma N.$$
 (Pick an open neighborhood in $S$ about each point in $S$ to obtain the above union.)  Thus $\Delta$ is closed if and only if it is of the form:
 $$\Delta = \bigcap_{N \lhd \Gamma} \bigcup_{\{\gamma \, | \, \gamma N \cap \Delta \ne \emptyset\}} \gamma N,$$
 When $\Delta$ is a subgroup,
 $$\Delta_N \defeq \bigcup_{\{\gamma \, | \, \gamma N \cap \Delta \ne \emptyset\}} \gamma N$$
 is an open subgroup for every fixed $N$. Indeed if
 $\gamma_{1} n_1, \gamma_{2}n_2 \in \Delta_N$ then $(\gamma_{1} n_1)^{-1} \gamma_{2}n_2 \in
 \gamma_1^{-1} \gamma_2 N$.  We can pick $n'_1, n'_2 \in N$ so that $\gamma_1 n'_1, \gamma_2 n'_2 \in \Delta$.  But $(\gamma_1 n'_1)^{-1} \gamma_2 n'_2 \in \gamma_1^{-1} \gamma_2 N$.  So the coset $\gamma_{1}^{-1} \gamma_{2} N$ must appear in the union defining
 $\Delta_N$. Finally the trivial coset $N$ appears in $\Delta_N$ in order to account for the
 identity element in $\Delta$.
 \end{remark}

 \begin{pf}
 Let $\Gamma < \PSL_2(\C)$ be a lattice, $\Delta < \Gamma$ a finitely generated subgroup. We
 first have to prove that the pro-normal topology is well defined on $\Gamma$. In search of a
 contradiction, assume that $N_1$ and $N_2$ are two non-trivial normal subgroups with $N_1 \cap
 N_2 = \trivgp$. By Borel's density theorem $\Gamma$ is Zariski dense in $\PSL_2(\C)$. Since
 $\PSL_2(\C)$ is simple, $N_1$ and $N_2$ are also Zariski dense. But $[N_1,N_2] < N_1 \cap N_2
 = \trivgp$, and this commutativity extends to the Zariski closure, proving that $\PSL_2(\C)$
 is commutative, which is absurd.

 When $\Delta$ is of finite index, it is clearly closed in the pro-normal topology and even in
 the pro-finite topology. Thus from here on we will assume that $\Delta$ is of infinite index
 and we will prove the second, stronger, assertion of Theorem \ref{thm:main}.
 We distinguish between two cases:

  \vskip 4pt {\noindent}\emph{Case (1):} $\Delta$ is geometrically finite. \vskip 2pt

  Recall we are assuming that $\Delta$ is of infinite index. In this case $\Delta$ cannot
  contain a non-trivial normal subgroup of $\Gamma$. This is because the limit set of a
  geometrically finite subgroup has measure zero, whereas $\Gamma$ and a non-trivial normal
  subgroup of $\Gamma$ have the same limit set, which is the whole sphere at infinity.

  We will prove that
  \begin{equation} \label{eqn:DefDelta}
  \Delta = \overline{\Delta} \ \defeq \ \bigcap_{ \{ N \lhd \Gamma, N \nleqslant \Delta \} } \Delta N.
  \end{equation}
  Assume by way of contradiction that $\overline{\Delta} \ne \Delta$ and let $\gamma \in
  \overline{\Delta} \setminus \Delta$. Applying Proposition \ref{main prop}, we find an infinite
  normal subgroup $N \lhd \Gamma$ such that $\Delta \cap N = \text{id}$, and therefore $\Delta N
  = \Delta \ltimes N$.  By the definition of $\overline{\Delta}$ we know that $\gamma \in
  \overline{\Delta} \le \Delta N$. Let $\gamma = \delta n$ be the unique way of factoring
  $\gamma$ into a product with $n \in N$ and $\delta \in \Delta$. Since $\Gamma$ is a finitely
  generated linear group it is residually finite and we can find a normal finite index subgroup
  $N' \lhd \Gamma$ such that $n \not \in N'$. Consider the group $M = N \cap N'$.  $M$ is
  nontrivial and intersects $\Delta$ trivially. Using the definition of $\overline{\Delta}$
  again, $\gamma \in \overline{\Delta} \le \Delta M$ so we can write $\gamma = \delta' m$ with
  $\delta' \in \Delta, m \in M$. But by construction $n \not \in M$, so $n \ne m$.  This
  contradicts the uniqueness of the factorization $\gamma = \delta n$.

   \vskip 4pt {\noindent}\emph{Case (2):} $\Delta$ is geometrically infinite. \vskip 2pt
   In this case we argue that $\Delta$ is pro-finitely closed in $\Gamma$, i.e. that it is the
   intersection of subgroups of finite index. This will prove the theorem because every finite
   index subgroup contains a normal subgroup (making it open in the pro-normal topology) and because $\Delta$, being of infinite index, cannot
   be equal to a finite index subgroup.

   Let us find groups $H$ and $F$ as in Lemma \ref{lemma (1)}.  Let $G = \langle H , \Delta
   \rangle$ be the group generated by $H$ and $\Delta$. Since $G$ is finite index in $\Gamma$ it
   is enough to show that $\Delta$ is closed in the pro-finite topology on $G$. Factoring out by
   the normal subgroup $F \lhd G$, it is enough to show that the finite subgroup $\Delta/F$ is
   pro-finitely closed in $G /F$. Note that $G/F$ is virtually cyclic and in particular it is
   residually finite.

   It remains only to observe that any finite subgroup is closed in the pro-finite topology on a
   residually finite group.  Let $H_i \lhd (G/F)$ be finite index normal subgroups with
   trivial intersection.  We claim that $ \Delta/F = \cap_i H_i (\Delta/F).$ Indeed let
    $g = h_1 \delta_1 = h_2 \delta_{2}  = \ldots = h_i \delta_{i}  = \ldots$ be any element of this
    intersection, where $h_i \in H_i$ and $\delta_i \in \Delta/F$.  But there are only finitely many ways to write this element in the form $h_i
    \delta_i$, because $\Delta/F$ is finite. In particular there exists some $h \in G/F$ such that
    $h_i = h$ for infinitely many $i'$s. But $\cap H_i = \trivgp$ so $h = \text{id}$ and $g \in
    \Delta/F$.
    \end{pf}
    
    \begin{remark} \label{rem:pro_char} One can also define the pro-characteristic topology as the topology generated by all the characteristic subgroups of $\Gamma$. This topology is coarser than the pro-normal topology but still finer than the pro-finite topology. We have in fact proved the stronger statement that every finitely generated subgroup $\Delta < \Gamma$ is closed in the pro-characteristic topology on $\Gamma$. Indeed if $\Delta$ is geometrically infinite or of finite index then it is closed even in the pro-finite topology. When $\Delta$ is a geometrically finite subgroup of infinite index, change the definition of the group $\overline{\Delta}$ so that the intersection in Equation (\ref{eqn:DefDelta}) is only over characteristic subgroups. The rest of the argument works just as before because, as indicated in Remark \ref{rem:characteristic}, the group $N$ can be chosen to be characteristic in $\Gamma$, and the normal subgroup $N'$ can be replaced by a subgroup that is characteristic but still of finite index.
    \end{remark}
    \bibliography{pronormal.biblio}
    \end{document}